\documentclass[10pt,
final,
conference,
twocolumn,
compsoc,
 letterpaper]{IEEEtran} 

\usepackage{srcltx}
\usepackage[psamsfonts]{amsfonts}
\usepackage{amsmath}
\usepackage{amssymb}
\usepackage{amsfonts}
\usepackage{graphics}
\usepackage{tikz} 
\usetikzlibrary{fit,positioning}
\usetikzlibrary{arrows,positioning,automata}
\usepackage{psfrag}
\usepackage{epsfig}
\usepackage{subfigure}
\usepackage{array}
\usepackage{algorithm}
\usepackage{algorithmic}
\usepackage{pifont}
\usepackage{srcltx}
\usepackage[psamsfonts]{amsfonts}
\usepackage{makeidx}  
\usepackage{bbm}
\usepackage{hhline}
\usepackage{eufrak}
\usepackage{yfonts}
\usepackage{color}
\usepackage{url}
\usepackage{dsfont}
\usepackage{caption}
\usepackage[square, comma, sort&compress, numbers]{natbib}
\usepackage{balance} 

\newtheorem{thm}{Theorem}

\newtheorem{lem}{Lemma}

\newtheorem{define}{Definition}

\IEEEoverridecommandlockouts

\interdisplaylinepenalty=2500 \algsetup{indent=2em}
\title{Robust Target Localization Based on Squared Range Iterative Reweighted Least Squares}
\author{Alireza Zaeemzadeh, Mohsen Joneidi, Behzad Shahrasbi, and Nazanin Rahnavard\\
 \thanks{This material is based upon work supported by the National Science Foundation under Grant No. ECCS-1418710 and Grant No. CCF-1718195.}
School of Electrical Engineering and Computer Science \\
University of Central Florida, Orlando, FL 32816, USA\\
Emails: {\{zaeemzadeh, joneidi, behzad, and nazanin\}@eecs.ucf.edu}}
\begin{document}
\renewcommand{\textfraction}{0}
\maketitle

\begin{abstract}
In this paper, the problem of target localization in the presence of outlying sensors is tackled. This problem is important in practice because in many real-world applications the sensors might report irrelevant data unintentionally or maliciously.
The problem is formulated by applying robust statistics techniques on squared range measurements and two different approaches to solve the problem are proposed. 
The first approach is computationally efficient; however, only the objective convergence is guaranteed theoretically. On the other hand, the whole-sequence convergence of the second approach is established. To enjoy the benefit of both approaches, they are integrated to develop a hybrid algorithm that offers computational efficiency and theoretical guarantees. \par

The algorithms are evaluated for different simulated and real-world scenarios. The numerical results show that the proposed methods meet the Cr\`amer-Rao lower bound (CRLB) for a sufficiently large number of measurements. When the number of the measurements is small, the proposed position estimator does not achieve CRLB though it still outperforms several existing localization methods. 

\end{abstract}
\begin{IEEEkeywords}
Target localization, robust localization, robust statistics, iterative reweighted least squares, generalized trust region subproblems
\end{IEEEkeywords}

\section{Introduction}\label{intro}
The problem of localization arises in different fields of study such as wireless networks, navigation, surveillance, and acoustics \cite{Ekim11Robust,Beck08range,Destino11onML}. There are many different approaches to localization based on various types of measurements such as range and squared-range (SR), time-of-arrival (ToA), time-difference-of-arrival (TDoA), two-way time-of-flight (TW-ToF), direction-of-arrival (DoA), and received-signal-strength (RSS)  \cite{Beck08range,Jiang07RobustTOA,
Jamali11SparsityTDOA,
Zhang09RobustTWTOF,Wang11NewApproach}. 

In \cite{Beck08range}, localization from range measurements and range-difference measurements are considered and least-squares (LS) estimators are exploited. Authors in \cite{Ekim11Robust,Beck08range,Destino11onML,Wang11NewApproach} have established methods to find the exact or approximate solution in the maximum likelihood (ML) framework. Usually finding the solution for ML  estimators is a difficult task or computationally burdensome \cite{Destino11onML,Wang11NewApproach}. \par

In this paper, the problem of \emph{robust target localization} is considered. 
In sensor networks, some nodes may report faulty data to the processing node unintentionally or maliciously. This may occur because of network failures, low battery, physical obstruction of the scene, and attackers. Thus, the processing node should not simply aggregate measurements from all sensors. It is more efficient to disregard the outlier measurements and localize the target based on reliable measurements.\par

There are different approaches toward robust localization. 
The method in \cite{Jiang07RobustTOA} is obtained by modeling the ToA estimation error as Cauchy-Lorentz distribution. In \cite{Champagne14DistLoc}, robust statistics, and specifically Huber norm, is exploited to localize sensors in a network in a  distributed manner using the location of a subset of nodes.  Authors in \cite{Zhang09RobustTWTOF} try to minimize the worst-case likelihood function and employ semidefinite relaxation to attain the estimate using TW-ToF measurements. 
The authors in \cite{Zoubir13RIN} have developed a robust geolocation method by estimating the probability density function (PDF) of the measurement error as a summation of Gaussian kernels. This method works best when the measurement error is drawn from a Gaussian mixture PDF. \par 

In this paper, the goal is to localize a single target in the presence of outlier range measurements in a centralized manner. We aim to achieve \emph{outlier distributional robustness}, which means the estimator performs well for different outlier probability distributions. A least squares methodology is applied to the squared range measurements. Although, this formulation is not optimal in the ML sense \cite{Destino11onML}, it provides us with the opportunity to find the estimate efficiently.


The contributions of this work can be summarized as follows. First, a robust optimization problem is formulated, which disregards unreliable measurements, using squared-range formulation. Next, two different algorithms are proposed to find the solution of the optimization problem.  
In the first algorithm, which is based on iteratively reweighted least squares (IRLS), the proposed optimization problem is transformed into a special class of optimization problems, namely Generalized Trust Region Subproblems (GTRS) \cite{More93GTRS}. Numerical simulations show that this algorithm has fast objective convergence. However the whole-sequence convergence is not established theoretically. \par 

The second algorithm is based on gradient descent. This algorithm is globally convergent, but needs more iterations to converge. By using these two algorithms, we proposed a hybrid method, which has desirable theoretical and practical features, such as fast whole-sequence convergence. \par 

The rest of this paper is organized in the following order. In Section~\ref{model}, the system model is introduced. Section~\ref{robustlocalization} describes the robust localization problem and two methods to tackle the problem are presented. 
Section~\ref{results} presents the simulation results and finally Section~\ref{conclusions} draws conclusions.\par

\section{System Model}\label{model}
Since the problem of source localization arises in different fields such as wireless networks, surveillance, navigation, and  acoustics, a general system model is exploited. In the generalized model, the system is comprised of $R$ sensors, with known locations, and the location of the target is estimated using the range measurements reported by these sensors. A central processing node collects the measurements and computes the location of the target.\par
Each sensor reports a range estimate, denoted by $r_i$, given by
\begin{equation}\label{measurements}
r_i = \| \boldsymbol{x} - \boldsymbol{a_i} \|_2 + v_i, \qquad  i = 1,...,R,
\end{equation}
where $\| . \|_2$ denotes Euclidean distance, $\boldsymbol{x} \in \mathbb{R}^n$ is the coordinates of the target, $\boldsymbol{a_i} \in \mathbb{R}^n$ is the location of the $i^{th}$ sensor and  $v_i$ models the measurement error. It is clear that for the aforementioned applications $n = 2$ or $3$. \par

The measurement errors $v_i$ are assumed to be independent and identically distributed random variables. To model the outlier measurements, a two-mode mixture PDF is assigned to the measurement errors, which can be written as:
\begin{equation}\label{noise-pdf}
p_V(v) = (1-\beta)\mathcal{N}(v;0,\sigma^2) + \beta \mathcal{H}(v).
\end{equation}

In other words, measurement errors are drawn from the distribution $\mathcal{N}(v;0,\sigma^2)$ with probability $1-\beta$ or the distribution  $\mathcal{H}(v)$ with probability $\beta$. $\mathcal{N}(v;0,\sigma^2)$ models the measurement noise for the outlier-free measurements, which is assumed to be a zero mean Gaussian distribution with variance $\sigma^2$, and $\mathcal{H}(v)$ models the outlier errors. Thus, the probability $\beta$ denotes the ratio of outlier measurements to all the measurements, also known as the \emph{contamination ratio}. The outlier error distribution, $\mathcal{H}(v)$, is commonly modeled with a Uniform distribution \cite{Hussain12Uniform,Nawaz11Uniform}, a shifted Gaussian distribution \cite{Gustafsson05positioning,Zoubir12nonparametric,Zoubir13RIN}, a Rayleigh distribution \cite{Zoubir12nonparametric}, or an exponential distribution \cite{Chen99exponential}. However, it is worthwhile to mention that our proposed method does not rely on the distribution of $\mathcal{H}(v)$. \par 




Here, the goal is to estimate $\boldsymbol{x}$ using the measurements $r_i \ i = 1, \ldots R$, while disregarding the measurements from outlier sensors. The processing node has no information about the number of the outlier sensors and the distribution of outlier measurements. Moreover, it is assumed that all the reported measurements including the noisy and irrelevant measurements are positive. For that, we exploit robust statistics and propose methods to obtain the solution. \par 

\section{Robust Localization From Squared Range Measurements}\label{robustlocalization}
In this section, a localization method is developed by applying robust statistics to the squared range measurements. Although this formulation is not optimal in the ML sense, unlike the methods based on range measurements, the solution can be attained easily. \par


The conventional square-range-based least squares (SR-LS) formulation is as follows \cite{Beck08range}:
\begin{equation}\label{SRLS}
\begin{aligned}
& \underset{\boldsymbol{x}}{\text{minimize}  }
& & \sum_{i = 1} ^{R}  ( \| \boldsymbol{x} - \boldsymbol{a_i} \|_2^2 - r_i^2)^2 .  \\
\end{aligned}
\end{equation}

It is clear that the problem stated in (\ref{SRLS}) is not convex. However, we can transform (\ref{SRLS}) into a special class of optimization problems by reformulating it as a constrained minimization problem given by \cite{Beck08range,More93GTRS}
\begin{equation}\label{SRLSconstrained}
\begin{aligned}
& \underset{\boldsymbol{x},\alpha}{\text{minimize}  }
& & \sum_{i = 1} ^{R}  (\alpha -2\boldsymbol{a_i}^T\boldsymbol{x} + \|\boldsymbol{a_i}\|^2 - r_i^2)^2 , \\
& \text{subject to}
& & \|\boldsymbol{x}\|^2 = \alpha.
\end{aligned}
\end{equation}

It is worthwhile to mention that $\alpha$ is also an outcome of the optimization procedure, not a parameter to be set. In this formulation, the unreliable measurements from outlier sensors affect the accuracy of localization significantly. 
We plan to use robust statistics to decrease the sensitivity of the estimator to the common assumptions. Here, \emph{robustness} signifies insensitivity to small deviation from the common assumption, which is the Gaussian distribution for noise. In (\ref{noise-pdf}), the parameter $\beta$ represents the deviation from this assumption. Our goal is to deal with the unknown distribution $\mathcal{H}(v)$ and to achieve \emph{distributional robustness}.\par 

As described in \cite{huber2011robust}, a proposed statistical procedure should have the following features. It must be \emph{efficient}, in the sense that it must have an optimal or near optimal performance at the assumed model, i.e., the Gaussian distribution for noise. It must be \emph{stable}, i.e., robust to small deviations from the assumed model. Also, in the case of \emph{breakdown}, or large deviation from the model, a catastrophe should not occur. In the numerical experiments, we will look for these features in the proposed methods. \par 

The general recipe to robustize any statistical procedure is to decompose the observations to fitted values and residuals \cite{huber2011robust}. In our proposed methods, we will try to find the residuals and re-fit iteratively until convergence is obtained. Each term of summation in (\ref{SRLSconstrained}) corresponds to the residual from a single sensor. These residuals can be exploited to re-fit the observations iteratively. \par 

Specifically, we use the residuals to assign weights to each observation. If an observation is fitted to the model, it should have a larger weight in the procedure of decision making. Inspired by \cite{Daubechies10IRLSconvergence}, we define the objective function as:

\begin{equation}
\label{AuxObj}
\mathcal{J}(\boldsymbol{y},\boldsymbol{w}) = \sum_{i = 1} ^R w_i(\boldsymbol{\tilde{a_i}}^T\boldsymbol{y} - b_i)^2 + \sum_{i = 1} ^R \epsilon^2 w_i - \ln{w_i},
\end{equation}
where $\boldsymbol{\tilde{a_i}}^T = \left[ \begin{array}{cc} -2\boldsymbol{a_i}^T & 1  \end{array} \right]$, $\boldsymbol{y} = \left[ \begin{array}{cc} \boldsymbol{x} & \alpha  \end{array} \right]^T$ , $b_i = r_i^2 - \|\boldsymbol{a_i}\|^2$, and $\boldsymbol{w} \in \mathbb{R}^R$ is the weight vector with $w_i > 0, \forall i$. The value of the parameter $\epsilon$ is a function of the standard deviation of the noise, we set $\epsilon = 1.34\sqrt[]{3}\sigma$ based on the discussion presented in Appendix \ref{app-eps}. 

The first summation of the objective function (\ref{AuxObj}) is the weighted version of the objective in (\ref{SRLSconstrained}). The other terms are added in such a way that result in the commonly used class of M-estimators known as Geman-McClure (GM) function \cite{Pennacchi08Mestimators,geman1987statistical}. The aim of GM function is reduce the effect of large errors, by interpolating between $\ell_2$ and $\ell_0$
norm minimizations \cite{Pennacchi08Mestimators}. There are other M-estimatiors with similar behavior as the Geman-McClure such as Tukey, Welsch, and Cauchy estimators. These types of M-estimators are known to be more robust to large errors than Huber M-estimator \cite{Pennacchi08Mestimators}. The desirable feature of Huber function is the convexity, unlike all the other mentioned estimators. However, our numerical results show that the proposed algorithms  perform well for different scenarios and different values of contamination ratio. \par 

Our goal is to minimize $\mathcal{J}(\boldsymbol{y},\boldsymbol{w} )$ over $\boldsymbol{y}$ and $\boldsymbol{w}$. Specifically, we are solving the following optimization problem:
\begin{equation}
\label{Eq-FullOpt}
\begin{aligned}
& 
& &  \underset{\boldsymbol{y},\boldsymbol{w}}{\text{minimize}  }
& & & \mathcal{J}(\boldsymbol{y},\boldsymbol{w} ), \\
& 
& & \text{subject to} 
& & & \boldsymbol{y}^T\boldsymbol{D}\boldsymbol{y} +2\boldsymbol{f}^T\boldsymbol{y} = 0, \\
& 
& & 
& & & w_i > 0, \forall i,
\end{aligned}
\end{equation}

where
\begin{equation}\label{SRL1Dfy}
\boldsymbol{D} = \left[ \begin{array}{cc} \boldsymbol{I}_n & \boldsymbol{0}_{n\times1} \\ \boldsymbol{0}_{1\times n} & 0 \end{array} \right],\, \boldsymbol{f} = \left[ \begin{array}{c}  \boldsymbol{0}_{n\times1}\\ -0.5 \end{array} \right].
\end{equation}

Our algorithms will exploit an alternative approach to update the weights and  $\boldsymbol{y}$. We initialize by taking $w_i^{(0)} = 1, \forall i$. Then at the $k^{th}$ iteration, the following optimization problem is solved to update the value of $\boldsymbol{y}$:

\begin{equation}
\label{Opt1}
\begin{aligned}
& \boldsymbol{y}^{(k+1)} = 
& & \arg \min 
\qquad \mathcal{J}(\boldsymbol{y},\boldsymbol{w}^{(k)} ), \\
& 
& & \text{subject to} 
\qquad \boldsymbol{y}^T\boldsymbol{D}\boldsymbol{y} +2\boldsymbol{f}^T\boldsymbol{y} = 0.
\end{aligned}
\end{equation}


Likewise, the weights are updated as follows:
\begin{equation}
\label{Opt2}
\begin{aligned}
& \boldsymbol{w}^{(k+1)} = 
& & \arg \min 
\qquad \mathcal{J}(\boldsymbol{y}^{(k+1)},\boldsymbol{w} ), \\
& 
& & \text{subject to} 
\qquad w_i > 0, \forall i.
\end{aligned}
\end{equation}



This problem is convex and the global minimizer can be obtained easily. As a result, the weights are given by:
\begin{equation}\label{SRIRLSW}
\begin{aligned}
& 
&& w_{i}^{(k)} = \frac{1}{ ({e_i^{(k)}})^2 + \epsilon^2},\\
& \text{where}
&& e_i^{(k)}= \boldsymbol{\tilde{a_i}}^T\boldsymbol{y}^{(k)} - b_i.
\end{aligned}
\end{equation}

Choosing such weights is common in iteratively reweighted least square (IRLS) methods \cite{Chartrand08IRLS,huber2011robust,Pennacchi08Mestimators,BoloursazMashhadi2017LevelThresholding,Zaeemzadeh2015}. \par 

In robust statistics terms, the measurements are decomposed into the fitted values $\boldsymbol{y}^{(k)}$ and residuals $\boldsymbol{e}^{(k)}$ at each iteration $k$. Then, the residuals are exploited to tune the weights of the observations. For large residuals, i.e., $e_i \gg \epsilon$, each term of the first summation in (\ref{AuxObj}), tends to $1$ . Similarly, for small residuals each term in summation tends to zero. In other words, we are minimizing the number of the observations with large residuals. \par

 
Now, two different approaches to find the solution of (\ref{Opt1}) are introduced. In the first approach, we show that (\ref{Opt1}) can be mapped into a special class of optimization problems known as Generalized Trust Region Subproblems (GTRS) \cite{More93GTRS}. Then at each iteration, the exact solution is derived by employing the GTRS formulations. In the second approach, a method based on gradient descent is introduced to solve the problem. This method is not as computationally efficient as the first approach, but offers an array of desirable theoretical features. \par 

\subsection{The Squared Range Iterative Reweighted Least Squares (SR-IRLS) Approach}\label{SR-IRLS}
The optimization problem in (\ref{Opt1}) can be formulated in the matrix form as:
\begin{equation}\label{SRL1matrix}
\begin{aligned}
& \underset{\boldsymbol{y}}{\text{minimize}  }
& & (\boldsymbol{A}\boldsymbol{y} - \boldsymbol{b})^T\boldsymbol{W}^{(k-1)}(\boldsymbol{A}\boldsymbol{y} - \boldsymbol{b}),   \\
& \text{subject to}
& & \boldsymbol{y}^T\boldsymbol{D}\boldsymbol{y} +2\boldsymbol{f}^T\boldsymbol{y} = 0,
\end{aligned}
\end{equation}
with
\begin{equation}\label{SRL1Ab}
\boldsymbol{A} = \left[ \begin{array}{cc} -2\boldsymbol{a_1}^T & 1 \\ \vdots & \vdots \\ -2\boldsymbol{a_R}^T & 1  \end{array} \right],\  \boldsymbol{b} = \left[ \begin{array}{c} r_1^2 - \|\boldsymbol{a_1}\|^2 \\ \vdots \\ r_R^2 - \|\boldsymbol{a_R}\|^2 \end{array} \right] ,
\end{equation}
and $\boldsymbol{W}^{(k)}$ is a diagonal weighting matrix in the $k^{th}$ iteration and $w_{i}^{(k)}$ is the $i^{th}$ diagonal entry of $\boldsymbol{W}^{(k)}, i = 1, \ldots, R$. \par 

Note that in (\ref{SRL1matrix}), a quadratic objective function is being minimized subject to a quadratic equality constraint. This special class of optimization problems is called Generalized Trust Region Subproblem (GTRS) \cite{More93GTRS}. The equality constraint makes this optimization problem non-convex. However, it is shown that the global solution of GTRS problems can be obtained efficiently \cite{More93GTRS,Beck08range}. \par  
\begin{thm}
\label{thm-global-existence}
\textit {
Let $q: \mathbb{R}^n \rightarrow \mathbb{R}$ and $c: \mathbb{R}^n \rightarrow \mathbb{R}$ be quadratics and assume $\{ \boldsymbol{x} \in \mathbb{R}^n : c(\boldsymbol{x}) = 0 \}$ is not empty. If
\begin{equation}
\label{eq-C-Q-relation}
\boldsymbol{v} \neq 0 , \boldsymbol{v}^T \boldsymbol{C} \boldsymbol{v} = 0 \Rightarrow \boldsymbol{v}^T \boldsymbol{Q} \boldsymbol{v} > 0,
\end{equation}
where
$$
\boldsymbol{Q} = \nabla^2 q,
\boldsymbol{C} = \nabla^2 c,
$$then the optimization problem $\min \{ q(\boldsymbol{x}) : c(\boldsymbol{x}) = 0\}$ has a global minimizer.
}
\end{thm}

\begin{thm}
\label{thm-global-conditions}
\textit {
Let $q: \mathbb{R}^n \rightarrow \mathbb{R}$ and $c: \mathbb{R}^n \rightarrow \mathbb{R}$ be quadratics and assume that  $ \min \{ c(\boldsymbol{x}) : \boldsymbol{x} \in \mathbb{R}^n \} < 0 < \max \{ c(\boldsymbol{x}) : \boldsymbol{x} \in \mathbb{R}^n \}$ and $\nabla^2 c \neq \boldsymbol{0}$. A vector $\boldsymbol{x}^*$ is a global minimizer of problem $\min \{ q(\boldsymbol{x}) : c(\boldsymbol{x}) = 0\}$ if and only if $c(\boldsymbol{x}^*) = 0$  and there is a multiplier $\lambda^* \in \mathbb{R}$ such that the Kuhn-Tucker condition
$$
\nabla q (\boldsymbol{x}^*) + \lambda^* \nabla c(\boldsymbol{x}^*) = \boldsymbol{0}
$$
is satisfied with
$$
\nabla^2 q(\boldsymbol{x}^*) + \lambda^* \nabla^2 c(\boldsymbol{x}^*)
$$
positive semidefinite.
}
\end{thm}

Specifically, using Theorem \ref{thm-global-existence} and the definitions of $\boldsymbol{A}$, $\boldsymbol{W}^{(k)}$, and $\boldsymbol{D}$, we can easily verify that (\ref{eq-C-Q-relation}) holds for the proposed optimization problem in (\ref{SRL1matrix}). Thus, the optimization problem (\ref{SRL1matrix}) has a global minimizer for all the iterations. Also by using Theorem \ref{thm-global-conditions}, $\boldsymbol{y}^{(k)}$ is an optimal solution of (\ref{SRL1matrix}) if and only if there exists $\lambda \in \mathbb{R}$ such that:
\begin{equation}\label{GTRS}
\begin{aligned}
& (\boldsymbol{A}^T \boldsymbol{W}^{(k-1)} \boldsymbol{A} + \lambda  \boldsymbol{D})\boldsymbol{y}^{(k)} = \boldsymbol{A}^T  \boldsymbol{W}^{(k-1)}  \boldsymbol{b} - \lambda  \boldsymbol{f},\\
& {\boldsymbol{y}^{(k)}}^T\boldsymbol{D}\boldsymbol{y}^{(k)} +2\boldsymbol{f}^T\boldsymbol{y}^{(k)} = 0, \\
& \boldsymbol{A}^T \boldsymbol{W}^{(k-1)} \boldsymbol{A} + \lambda  \boldsymbol{D} \succeq 0. \\
\end{aligned}
\end{equation}

The last expression means that $\boldsymbol{A}^T \boldsymbol{W}^{(k-1)} \boldsymbol{A} + \lambda  \boldsymbol{D}$ is positive semidefinite. 
The first two equalities in (\ref{GTRS}) can be exploited to obtain a solution for $\lambda$, i.e. $\lambda^*$. 
To ensure that $\boldsymbol{A}^T \boldsymbol{W}^{(k-1)} \boldsymbol{A} + \lambda^*  \boldsymbol{D}$  is positive semidefinite, it is easy to show that we need to seek for $\lambda^*$ in the interval
\begin{equation}\label{Interval}
\lambda^* \geq - \frac{1}{\lambda_1(\boldsymbol{D},\boldsymbol{A}^T \boldsymbol{W}^{(k-1)} \boldsymbol{A})} ,
 \end{equation}
where $\lambda_1(\boldsymbol{D},\boldsymbol{A}^T \boldsymbol{W}^{(k-1)} \boldsymbol{A})$ is the largest generalized eigenvalue of the matrix pair $(\boldsymbol{D},\boldsymbol{A}^T \boldsymbol{W}^{(k-1)} \boldsymbol{A})$. 
It is shown that if (\ref{eq-C-Q-relation}) holds, then $\boldsymbol{A}^T \boldsymbol{W}^{(k-1)} \boldsymbol{A} + \lambda  \boldsymbol{D} \succeq 0$ for some $\lambda \in \mathbb{R}$ \cite[Theorem 2.2]{More93GTRS}. Moreover, the resulting characteristic function needed to be solved to find $\lambda^*$ is strictly decreasing over this interval \cite[Theorem 5.2]{More93GTRS}. Thus, at each iteration, $\lambda^*$ can be obtained using a bisection algorithm. The interval for starting point of the bisection algorithm is specified as $(\lambda_l,\infty)$, where $\lambda_l = \max \{- (\boldsymbol{A}^T \boldsymbol{W}^{(k-1)} \boldsymbol{A})_{ii}, i = 1, \ldots, n\}$ \cite{More93GTRS}. \par 

Then, $\boldsymbol{y}$ is updated using the estimated $\lambda^*$. Algorithm \ref{SRIRLSalg} illustrates the procedure to calculate the estimate of (\ref{SRL1matrix}) using the equations in (\ref{GTRS}). The convergence of the algorithm is analyzed in Theorem \ref{convergencetheorem}.
\begin{algorithm}
\caption{Calculating the SR-IRLS estimate}
\label{SRIRLSalg}
\algsetup{
linenosize=\small,
linenodelimiter=:
}
\begin{algorithmic}[1]

\REQUIRE $\boldsymbol{a_i}$, $\boldsymbol{r_i}  $ for $i = 1, \dots , R$, $\epsilon$, maximum number of iterations $maxIter$, and the convergence tolerance $\Delta$. 
\STATE $ \textbf{Compute } \boldsymbol{A},\boldsymbol{b},\boldsymbol{D}$, and $\boldsymbol{f}$ using (\ref{SRL1Ab}) and (\ref{SRL1Dfy}).
\STATE $\textbf{Initialize } w_i^{(0)} = 1, \forall i,$ and $ k = 1$.
\REPEAT
\STATE $\lambda_l = \max \{- (\boldsymbol{A}^T \boldsymbol{W}^{(k-1)} \boldsymbol{A})_{ii}, i = 1, \ldots, n\}.$
\STATE $\indent  \textbf{Find } \lambda^*: $ solve $\small \boldsymbol{y}(\lambda)^{T}\boldsymbol{D} \boldsymbol{y}(\lambda) +2\boldsymbol{f}^T\boldsymbol{y}(\lambda) = 0$ using a bisection algorithm in interval $(\lambda_l,\infty)$, where $\small \boldsymbol{y}(\lambda) = (\boldsymbol{A}^T \boldsymbol{W}^{(k-1)} \boldsymbol{A} + \lambda  \boldsymbol{D})^{-1} (\boldsymbol{A}^T  \boldsymbol{W}^{(k-1)}  \boldsymbol{b} - \lambda  \boldsymbol{f})$.
\STATE $ \textbf{Update }\small \boldsymbol{y}: \boldsymbol{y}^{(k)} = \boldsymbol{y}(\lambda^*) $.
\STATE $  \textbf{Update } {\boldsymbol{w}^{(k)}}$ using (\ref{SRIRLSW}).
\UNTIL {Convergence, i.e., if $| \mathcal{J}(\boldsymbol{y}^{(k)},\boldsymbol{w}^{(k)}) - \mathcal{J}(\boldsymbol{y}^{(k-1)},\boldsymbol{w}^{(k-1)}) | < \Delta$ or the maximum number of iterations $maxIter$ is reached. }
 \end{algorithmic}
\end{algorithm}

\begin{thm}
\label{convergencetheorem}
\textit { The sequence $ \{ \mathcal{J}(\boldsymbol{y}^{(k)},\boldsymbol{w}^{(k)}) \} $ generated by by Algorithm \ref{SRIRLSalg} converges to a constant value and every limit point of the iterates $\{ \boldsymbol{y}^{(k)},\boldsymbol{w}^{(k)} \}$ is a stationary point of (\ref{Eq-FullOpt}).
}
\end{thm}

\begin{IEEEproof}
See Appendix \ref{conv-analysis}.
\end{IEEEproof}

Inspection of the algorithm reveals that the matrix inversions are only needed for $(n+1) \times (n+1)$ matrices, where $n$ is the space dimension and is equal to $2$ or $3$. Thus the main computational burden of the algorithm stems from the matrix multiplications. The per iteration complexity of the algorithm is $\mathcal{O}(R^2)$. Similarly, the growth rate for the legacy least square problem is also $\mathcal{O}(R^2)$. Thus, the main computational burden of the SR-IRLS algorithm arises from the number of the iterations. \par

Our numerical experiments show that the SR-IRLS method needs a few iterations to solve the problem. The convergence of the objective is also proven in Appendix \ref{conv-analysis}. However, due to the lack of convexity, the standard convergence analysis tools cannot be used to show the convergence of the whole-sequence of the iterates. The problem becomes more difficult when the objective function is not a linear or a quadratic function of the previous iterates. Thus, in Appendix \ref{conv-analysis}, the convergence of a \emph{subsequence} of the iterates to a critical point is proved, although the whole-sequence convergence is almost always observed.\par 

This motivates us to propose a globally convergent algorithm. In Section \ref{Sec-SR-GD}, an algorithm, referred to as SR-GD, is introduced to find the solution of (\ref{SRL1matrix}) based on gradient descent. Then, we will integrate SR-IRLS and SR-GD to derive a \emph{computationally efficient} and \emph{globally convergent} algorithm. \par 

\subsection{The Squared Range Gradient Descent (SR-GD) Approach}\label{Sec-SR-GD} 

In this section, a new algorithm for solving the optimization problem in (\ref{Opt1}) is proposed based on gradient descent (SR-GD), for which the convergence of the whole-sequence of the iterates has been proven theoretically \cite{xu2014globally}. For that, the Lipschitz continuity of the gradient of the objective function as well as the special form of the objective and the constraint are employed. The numerical experiments show that this algorithm needs more iterations to converge than the SR-IRLS. Our goal will be to employ SR-GD and SR-IRLS to propose a hybrid fast converging algorithm. \par 

Inspired by \cite{xu2014globally}, at each iteration, the value of $\boldsymbol{y}^{(k)}$ is updated as follows:

\begin{equation}\label{eq-general-bcd}
\begin{aligned}
& \boldsymbol{y}^{(k)} = \arg\min_{\boldsymbol{y}} 
& & \langle \nabla_{\boldsymbol{y}} \mathcal{J}(\hat{\boldsymbol{y}}^{(k)},\boldsymbol{w}^{(k-1)}) , \boldsymbol{y} - \hat{\boldsymbol{y}}^{(k)} \rangle \\
&
& & + l^{(k)} \| \boldsymbol{y} - \hat{\boldsymbol{y}}^{(k)} \|_2^2,   \\
& \text{subject to}
& & \boldsymbol{y}^T\boldsymbol{D}\boldsymbol{y} +2\boldsymbol{f}^T\boldsymbol{y} = 0,
\end{aligned}
\end{equation}
where $$\hat{\boldsymbol{y}}^{(k)} = \boldsymbol{y}^{(k-1)} + \omega^{(k)}(\boldsymbol{y}^{(k-1)} - \boldsymbol{y}^{(k-2)}), $$ and $l^{(k)}$ is the Lipschitz constant of $\nabla_{\boldsymbol{y}} \mathcal{J}(\boldsymbol{y},\boldsymbol{w}^{(k-1)})$ at the $k^{\text{th}}$ iteration. By the definition of Lipschitz continuity, we have $$ \| \nabla_{\boldsymbol{y}} \mathcal{J}(\boldsymbol{u},\boldsymbol{w}^{(k-1)}) - \nabla_{\boldsymbol{y}} \mathcal{J}(\boldsymbol{v},\boldsymbol{w}^{(k-1)}) \| \leq l^{(k)} \|\boldsymbol{u} - \boldsymbol{v} \|. $$\par

Intuitively, the first term of the objective finds the steepest descent, while the second term prevents large changes in the magnitude of the gradient. The Lipschitz constant of the gradient function limits the step size of the algorithm and the new estimate $\boldsymbol{y}^{(k)}$ is enforced to be around the prediction $\hat{\boldsymbol{y}}^{(k)}$. The prediction is constructed using the previous iterates and an extrapolation factor $\omega^{(k)} = \frac{1}{12} \sqrt{\frac{l^{(k-1)}}{l^{(k)}}}$ \cite{xu2014globally}. The update rule for $\boldsymbol{w}$ is the same as (\ref{SRIRLSW}). \par 

This problem is not convex as well, but authors in \cite{xu2014globally} have proven the convergence of the whole-sequence of the algorithm by exploiting the properties of the objective. \par 

It is easy to notice that the minimization problem stated in (\ref{eq-general-bcd}) is a GTRS problem. This is because a quadratic objective is minimized subject to a quadratic equality constraint. By exploiting the definition of $\boldsymbol{D}$ and $l^{(k)}$, we can show that (\ref{eq-C-Q-relation}) holds. Thus the optimization problem in (\ref{eq-general-bcd}) has global minimizer for all iterations. Also Theorem \ref{thm-global-conditions} states that $\boldsymbol{y}^{(k)}$ is an optimal solution of (\ref{eq-general-bcd}) if and only if there exists $\lambda \in \mathbb{R}$ such that:

\begin{equation}\label{GTRS-bcd}
\begin{aligned}
&
\begin{aligned}
& (l^{(k)} \boldsymbol{I}_{n+1} + \lambda  \boldsymbol{D})\boldsymbol{y}^{(k)} =
&& -\boldsymbol{A}^T  \boldsymbol{W}^{(k-1)}  ( \boldsymbol{A} \hat{\boldsymbol{y}}^{(k)} - \boldsymbol{b} ) \\
&
&& + l^{(k)}  \hat{\boldsymbol{y}}^{(k)} - \lambda  \boldsymbol{f},\\
\end{aligned}
&& \\
 & {\boldsymbol{y}^{(k)}}^T\boldsymbol{D}\boldsymbol{y}^{(k)} +2\boldsymbol{f}^T\boldsymbol{y}^{(k)} = 0, 
 && \\
& \lambda \geq \max \{ -l^{(k)}, - \frac{1}{\lambda_1(\boldsymbol{D},\boldsymbol{A}^T \boldsymbol{W}^{(k-1)} \boldsymbol{A})} \} 
\end{aligned}
\end{equation}

At each iteration, after finding the predicted value for the iterate $\hat{\boldsymbol{y}}^{(k)}$, the equality expressions in (\ref{GTRS-bcd}) are used to find $\lambda$ and to update the values of $\boldsymbol{y}$ and $\boldsymbol{w}$. We should look for the solution of $\lambda$ in an interval that satisfies the positive semidefiniteness constraint. Since (\ref{eq-C-Q-relation}) holds, this interval exists and the characteristic function is strictly decreasing over this interval \cite[Theorem 2.2, Theorem 5.2]{More93GTRS}. Algorithm \ref{alg-SRGD} shows the steps to find the solution of the localization problem using the SR-GD method.


\begin{algorithm}
\caption{Calculating the SR-GD estimate}
\label{alg-SRGD}
\algsetup{
linenosize=\small,
linenodelimiter=:
}
\begin{algorithmic}[1]

\REQUIRE $\boldsymbol{a_i}$, $\boldsymbol{r_i}  $ for $i = 1, \dots , R$, $\epsilon$, maximum number of iterations $maxIter$, and the convergence tolerance $\Delta$. 
\STATE $ \textbf{Compute } \boldsymbol{A},\boldsymbol{b},\boldsymbol{D}$, and $\boldsymbol{f}$ using (\ref{SRL1Ab}) and (\ref{SRL1Dfy}).
\STATE $\textbf{Initialize } \boldsymbol{W}^{(0)}$ with identity matrix, $\boldsymbol{y}^{(-1)} = \boldsymbol{y}^{(0)} = \boldsymbol{A}^\dag\boldsymbol{b}$, $l^{(0)} = 0$, and $ k = 1$.
\REPEAT
\STATE $\indent l^{(k)} = 2 \| \boldsymbol{A}^T \boldsymbol{W}^{(k-1)} \boldsymbol{A}\|_F $.
\STATE $\indent \omega^{(k)} = \frac{1}{12} \sqrt{\frac{l^{(k-1)}}{l^{(k)}}}$.
\STATE $\indent \hat{\boldsymbol{y}}^{(k)} = \boldsymbol{y}^{(k-1)} + \omega^{(k)}(\boldsymbol{y}^{(k-1)} - \boldsymbol{y}^{(k-2)})$.
\STATE $\indent  \textbf{Find } \lambda^*: $ solve $\small \boldsymbol{y}(\lambda)^{T}\boldsymbol{D} \boldsymbol{y}(\lambda) +2\boldsymbol{f}^T\boldsymbol{y}(\lambda) = 0$ using a bisection algorithm in interval $(-l^{(k)},\infty)$, where $\small \boldsymbol{y}(\lambda) = (l^{(k)} \boldsymbol{I}_{n+1} + \lambda  \boldsymbol{D})^{-1} (-\boldsymbol{A}^T  \boldsymbol{W}^{(k-1)}  ( \boldsymbol{A} \hat{\boldsymbol{y}}^{(k)} - \boldsymbol{b} ) + l^{(k)}  \hat{\boldsymbol{y}}^{(k)} - \lambda  \boldsymbol{f})$.
\STATE $ \textbf{Update }\small \boldsymbol{y}: \boldsymbol{y}^{(k)} = \boldsymbol{y}(\lambda^*). $
\STATE $  \textbf{Update } {\boldsymbol{w}^{(k)}}$ using (\ref{SRIRLSW}).
\UNTIL {Convergence, i.e., if $\| \boldsymbol{y}^{(k)} - \boldsymbol{y}^{(k-1)} \| < \Delta$ or the maximum number of iterations $maxIter$ is reached. }
 \end{algorithmic}
\end{algorithm}

The numerical experiments show that the SR-GD method needs more time to find the solution than SR-IRLS. This is due to the fact that in SR-GD, the value of the new iterate is bounded to be around the previous iterate, unlike the SR-IRLS method. \par 

To take advantage of the fast convergence of the SR-IRLS and the whole sequence convergence of the SR-GD, we propose a hybrid method. Specifically, we can start with the SR-IRLS method and update the iterates by steps stated in Algorithm \ref{SRIRLSalg}. After convergence of the objective function, which is proven in Appendix \ref{conv-analysis}, the update rules in Algorithm \ref{alg-SRGD} are employed to find the final solution. The performance, convergence rate, and computational cost of this hybrid method is evaluated and compared with other methods in Section \ref{results}. \par

\section{Numerical Results}\label{results}
In this section, we present the simulation results to evaluate the performance of our proposed methods. We will seek for the main features of a robust estimator, which are discussed in Section \ref{SRLS}. We will examine the performance of the algorithms at the assumed model ($\beta = 0$), small deviations from model (small $\beta$), and large deviations from the model (large $\beta$). Moreover, we check distributional robustness of the proposed algorithms, which means that the performance of the methods will be evaluated for different outlier noise distributions $\mathcal{H}(v)$. \par 

Two different simulation scenarios will be investigated. In Scenario I, a general system model is considered and the outlier measurements obey a uniform distribution, which models a harsh environment. In Scenario II, localization of a target in a cellular radio network is investigated. The geometry of sensors is taken from an operating network and the measurement errors are drawn from a Gaussian mixture distribution to model the non-line-of-sight (NLOS) measurements. \par 

The performances of the proposed methods are compared with existing least-square-based \cite{Beck08range} and robust \cite{Zoubir13RIN} methods. \par

\subsection{Scenario I} \label{subsec-scene1}
The simulation parameters are as follows, unless otherwise is stated. In a $4000 \times 4000$ $m^2$ area, there exist $10$ sensors trying to localize a target. The sensors and the target are distributed uniformly at random, The range measurements are corrupted by the additive white Gaussian noise with standard deviation of $\sigma = 55 \ m$. Moreover, among the sensors, there exist $4$ outlier sensors. The noise of the outlier sensors are uniformly distributed in range $[-4000\sqrt{2},4000\sqrt{2}]$. Mathematically speaking, the distribution of the measurement error is as follows:\par
\begin{equation}\label{eq-uniform-pdf}
p_V(v) = (1-\beta)\mathcal{N}(v;0,\sigma^2) + \beta \mathcal{U}(v;-D_{max},D_{max}),
\end{equation}
where $\mathcal{U}(v;-D_{max},D_{max})$ is a uniform distribution with support $[-D_{max},D_{max}]$, which is modeling the outlier measurements. $\mathcal{N}(v;0,\sigma^2)$ is a zero mean Gaussian distribution with variance $\sigma^2$. \par 

To ensure that all the range measurements are positive, we set the non-positive values equal to a small value, i.e. $10^{-5}$. Localization is performed in a $2$-dimensional space, i.e. $n = 2$. \par

The performances of the proposed methods SR-IRLS, SR-GD, and the hybrid version are compared with the performance of SR-LS \cite{Beck08range}, a least-square-based method, as well as a robust method, namely Robust Iterative Nonparametric (RIN) \cite{Zoubir13RIN}. The performances are compared according to the root mean square error (RMSE),
\begin{equation}\label{RMSE}
\sqrt{\frac{1}{n} \| \boldsymbol{x} - \widehat{\boldsymbol{x}}\|_2^2},
\end{equation}
averaged over sufficiently large random simulations. $\widehat{\boldsymbol{x}}$ is the estimated value of the target location $\boldsymbol{x}$.\par

\begin{figure}
\centering
\includegraphics[width=3in,angle=0]{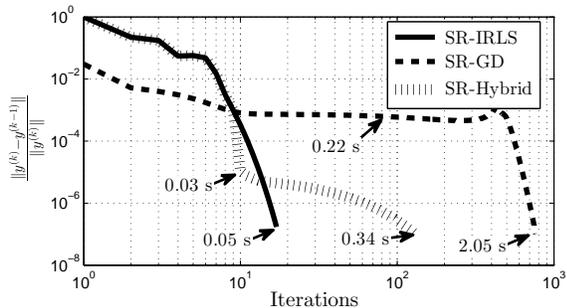}
\vspace{-1.5mm}
\caption{\small{Convergence of SR-IRLS, SR-GD, and the hybrid method. Labels show the execution time of different algorithms at different iterations. }}
\label{fig-convergence}
\vspace{-5mm}
\end{figure}

In our first numerical experiment, the convergence of SR-IRLS and SR-GD are compared. Figure \ref{fig-convergence} depicts $\frac{\| y^{(k)}  - y^{(k-1)}\|}{\| y^{(k)} \|}$ at different iterations. Moreover, the labels show the elapsed time for some of the iterations. Although the convergence of the SR-GD method is theoretically provable, Figure \ref{fig-convergence} shows that it needs more iterations and more time to converge. 
The hybrid version of the algorithm (SR-Hybrid) uses the update rules of SR-IRLS until the convergence of the objective function, then it employs the update rules of SR-GD. As a result, it needs less iterations than SR-GD, while its convergence is still theoretically provable. \par  

\begin{figure}
\centering
\includegraphics[width=3in,angle=0]{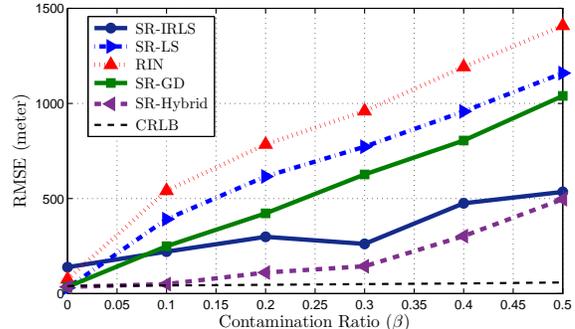}\vspace{-1.5mm}
\caption{\small{Robustness against outliers for $200$ Monte Carlo trials, $\sigma = 55$ and $R = 10$. Number of outlier sensors is set to $\beta \times R$. }}
\label{RMSEvsNumOutliers}
\vspace{-5mm}
\end{figure}

To study the influence of the number of outlier sensors, Figure \ref{RMSEvsNumOutliers} exhibits the RMSE of the estimate for different number of outlier sensors, or equivalently different values of $\beta$. In this study, the results are based on $200$ Monte Carlo (MC) trials. It is clear that as the number of outliers increases, the performance of the SR-LS method deteriorates significantly.
SR-IRLS and SR-GD perform closely for small values of $\beta$, but the difference becomes more noticeable as $\beta$ increases. This was expected since SR-GD is more likely to result in local optimum solutions caused by the outliers, because of the smooth convergence of the iterates. However, the hybrid version, which only uses small step size when it is sufficiently close to the limit point, performs the best for different values of contamination ratio. This figure shows that the proposed methods are efficient at the assumed method ($\beta = 0$) and stable for small deviations. Also, for large deviations, a catastrophe is not occurred. \par

To estimate the target location, the RIN method \cite{Zoubir13RIN} approximates the PDF of the measurement error with a summation of Gaussian kernels. For that, it needs a considerable number of measurements. Hence, unlike our proposed methods, it cannot produce proper results with $ R = 10$ measurements. Further, since the RIN method employs Gaussian kernels, it works most accurately when the measurement errors are drawn from a Gaussian mixture distribution (see Section \ref{subsec-scene2}). Using the Gaussian kernels decreases the distributional robustness of the RIN method significantly. \par 

\begin{figure}[h]
\centering     
\subfigure[Bias]{\label{subfig-biasvssensors}\includegraphics[width=3.1in]{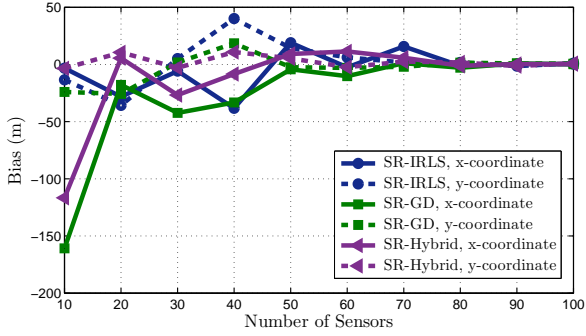}
\vspace{-5mm}}
\subfigure[RMSE]{\label{subfig-RMSEvssensors}
\includegraphics[width=3in]{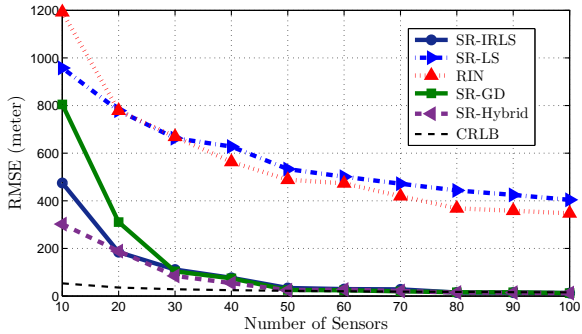}}
\subfigure[Running time]{\label{subfig-timevssensors}\includegraphics[width=3in]{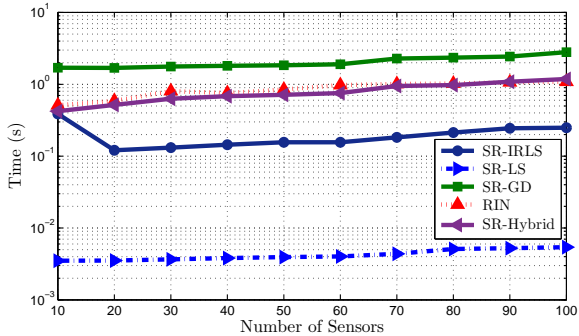}}
\caption{Performance of the localization methods versus number of sensors for $1000$ Monte Carlo trials and $\beta = 0.4$. }
\label{fig-vsNumofSensors}
\vspace{-5mm}
\end{figure}

To elaborate the point, Figure \ref{fig-vsNumofSensors} illustrates the impact of the number of sensors on performance of different methods. In this experiment, $40\%$ of the sensors are reporting unreliable data to the processing node, i.e. $\beta = 0.4$. This figure exhibits that the accuracy of the localization methods improves as the number of sensors increases. As it was expected, the RMSE of the estimates produced by the RIN method decreases significantly as the number of sensors increases. \par

Moreover, It is clear that the proposed methods meet the Cr\`amer-Rao lower bound (CRLB) for large number of measurements. From Figure \ref{subfig-biasvssensors} and Figure \ref{subfig-RMSEvssensors}, we can infer that the proposed methods are efficient for this simulation parameters, because they meet the CRLB and they are unbiased.  The CRLB is approximated by using Monte Carlo integration techniques explained in \cite{Zoubir13RIN}.  \par

Figure \ref{subfig-timevssensors} shows the running times\footnote{Running time reflects the time required to execute all the steps of the algorithms, including initialization, preprocessing, convergence, and post-processing. All simulations have been performed under MATLAB 2014a environment on a PC equipped with Intel Xeon E5-1650 processor (3.20 GHz) and 8 GB of RAM. } for different number of sensors. Clearly, the iterative methods requires more computation time than the least square method. Also, as it was expected and can be noticed in Figure \ref{fig-convergence}, the running time of the hybrid method is less than SR-GD, but more than SR-IRLS. \par

\begin{figure}
\centering
\includegraphics[width=3in,angle=0]{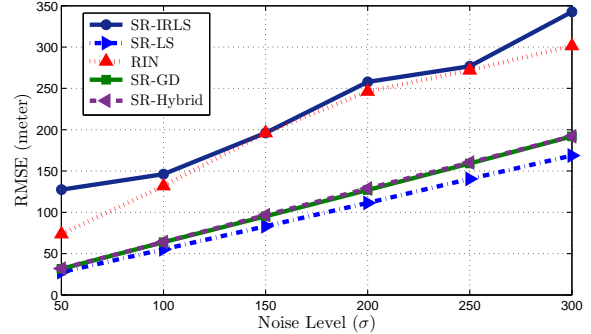}\vspace{-1.5mm}
\caption{\small{Comparison of the RMSEs in an environment with no outlier sensor, $\beta = 0$. }}
\label{fig-RMSEvsNoiseNoOutlier}
\vspace{-5mm}
\end{figure}

It is also worthwhile to compare the performance of the localization methods for the case when no sensor is reporting unreliable measurements and the range measurements are corrupted only by an additive Gaussian noise, i.e. $\beta = 0$. As it can be seen in Figure \ref{fig-RMSEvsNoiseNoOutlier}, the LS method outperform the robust methods. This was expected since the LS methods are particularly tailored to deal with Gaussian noise, while the robust methods are customized to handle the unreliable measurements. We are sacrificing efficiency for $\beta = 0$, to achieve stability in deviation from the model. However, it is easy to notice that the RMSE of the proposed robust methods is close to the RMSE of the SR-LS method, which implies the near optimal performance for Gaussian noise. \par

\subsection{Scenario II} \label{subsec-scene2}
\begin{figure}
\centering
\includegraphics[width=2.5in,angle=0]{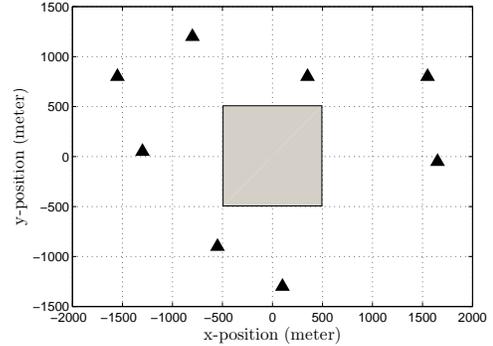}\vspace{-1.5mm}
\caption{\small{Geometry of the sensors, marked as triangle, and the city center area, marked as gray square, in a real world operating cellular radio network. }}
\label{fig-citycenter}
\vspace{-5mm}
\end{figure}

In this section, the problem of localizing a target in a radio cellular network is considered. The network consists $R = 8$ base stations (BSs), which are trying to estimate the location of a target in a city center area. The configuration of the BSs and the city center, as depicted in Figure \ref{fig-citycenter}, is taken from a realistic network \cite{Zoubir13RIN}. \par 

The outlier-free measurements are result of line-of-sight (LOS) sensings. On the other hand NLOS sensings produce unreliable measurements. Field trials have indicated that the measurement errors in harsh LOS/NLOS environments can be modeled as a 
Gaussian mixture distribution \cite{Zoubir13RIN},
\begin{equation}\label{eq-mixedgauss-pdf}
p_V(v) = (1-\beta)\mathcal{N}(v;0,\sigma^2) + \beta \mathcal{N}(v;\mu_{NL},\sigma_{NL}^2),
\end{equation}
where $\mathcal{N}(v;\mu_{NL},\sigma_{NL}^2)$ is a Gaussian distribution with mean $\mu_{NL}$ and variance $\sigma_{NL}^2$, modeling the NLOS measurements. \par 

For each BS, we obtain $K$ measurements and stack them in the measurement vector as follows:
\begin{equation}\label{eq-Kmeas-Ab}
\scriptsize
\boldsymbol{b} = \left[ \begin{array}{c} r_1(1)^2 - \|\boldsymbol{a_1}\|^2 \\ \vdots \\ r_1(K)^2 - \|\boldsymbol{a_1}\|^2 \\ \vdots \\ r_R(1)^2 - \|\boldsymbol{a_R}\|^2 \\ \vdots \\ r_R(K)^2 - \|\boldsymbol{a_R}\|^2 \end{array} \right] , \ 
\boldsymbol{A} = \left[ \begin{array}{cc} -2\boldsymbol{a_1}^T & 1 \\ \vdots & \vdots \\
-2\boldsymbol{a_1}^T & 1 \\ \vdots & \vdots \\ -2\boldsymbol{a_R}^T & 1  \\ \vdots & \vdots \\ -2\boldsymbol{a_R}^T & 1 \end{array} \right]. 
\end{equation}

In the simulations, it is assumed that each BS reports $ K = 20$ samples. The measurement errors are drawn from the distribution in (\ref{eq-mixedgauss-pdf}) with $\sigma = 55$, $\mu_{NL} = 380$, and  $\sigma_{NL} = 120$. The position of the target is uniformly generated in the city center area.\par 

\begin{figure}
\centering
\includegraphics[width=3in,angle=0]{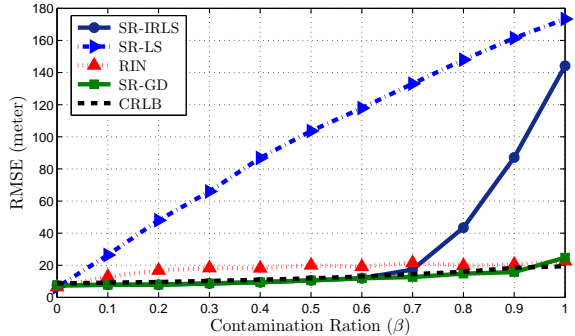}\vspace{-1.5mm}
\caption{\small{ Mean RMSE of different localization methods versus contamination ratio, for $100$ MC trials. }}
\label{fig-vsOutlierCityCenterScnario}
\vspace{-5mm}
\end{figure}

Figure \ref{fig-vsOutlierCityCenterScnario} illustrates the performance of different localization methods versus the contamination ratio for $0 \leq \beta  \leq 1$. This figure shows that SR-GD outperforms its competitors. Moreover, the hybrid version and SR-IRLS perform the same in this configuration and are able to handle NLOS measurements up to a certain amount and meet the CRLB up to a certain $\beta$. For, Large values of $\beta$, the SR-IRLS method breaks down, but still works better than the least square method. However, in this scenario SR-GD is able to localize the target for even large contamination ratios.
\par 

Moreover, the RIN method performs accurately in this scenario, in comparison with the previous scenario. In this scenario, the RIN can estimate the PDF of the error more accurately. This was expected because, firstly, we are collecting $R \times K = 160$ measurements, secondly, the measurement error has a mixture Gaussian distribution. As a result, the RIN method can produce a better estimate of the target location. With $160$ measurements, RIN is able to approximate the measurement error distribution. Thus, the its RMSE does not change considerably for different values of $\beta$. This fact is vividly clear for the extreme case. For $\beta = 1$, the RIN method is able to approximate $\mathcal{N}(v;\mu_{NL},\sigma_{NL}^2)$ as the PDF of the measurement error. As a result, this method  outperforms the competitors for the special case of $\beta = 1$. \par

\vspace{-2.5mm}
\section{Conclusions}\label{conclusions}
In this paper, we have considered the problem of localizing a single target in the presence of unreliable measurements with unknown probability distribution. For that, the squared-range formulation is exploited. To disregard the outlier measurements and find the estimate using the reliable measurements, we have used robust statistics. Then the problem is converted into a known class of optimization problems, namely GTRS, using the concepts in robust statistics. Two algorithms and a hybrid method are proposed to solve the problem. Convergence of the algorithms is analyzed theoretically. \par 

The simulation results suggest that the proposed methods outperform the existing methods, while providing a near optimal performance for Gaussian noise.\par


\appendices 
\section{Choosing the Parameter $\epsilon$}
\label{app-eps}

Here, we establish a connection between the objective function introduced in (\ref{AuxObj}) with Huber norm and use the results in robust statistics to tune $\epsilon$. The first summation in (\ref{AuxObj}) can be seen as an iterative approximation of $\sum_i\rho_\epsilon(e_i)$, where
\begin{equation}\label{Rho}
\rho_\epsilon(x)  = \frac{x^2}{x^2 + \epsilon^2}.
\end{equation}
$\rho_\epsilon(.)$ indicates a measurement as an outlier if the residual is greater than a threshold and this threshold is a function of $\epsilon$. Robustness to noise is improved by increasing the value of $\epsilon$, at the expense of losing robustness to the outlier measurements. Hence, as the variance of noise increases, we should assign a larger $\epsilon$ to $\rho_\epsilon(.)$.  To set the value of $\epsilon$, a link between the proposed problem and the Huber norm is established. \par 


In robust statistics \cite{huber2011robust}, Huber norm, $\rho_\tau^H(.)$, is utilized to disregard the outlier measurements. $\rho_\tau^H(.)$ is defined as
\begin{equation}
\scriptsize
\label{Huber-norm}
\rho_\tau^H(x) = \left\{ \begin{array}{lr} \frac{1}{2}x^2   & : | x | < \tau \\[2ex] \tau |x| - \frac{\tau^2}{2}  & : | x | \geq \tau \end{array} \right. 
\end{equation}

Assuming that the additive measurement noise is Gaussian, the estimator would be $95\%$ asymptotically efficient, meets Cr\`amer-Rao bound, by setting the parameter $\tau$ to $1.34\sigma$, where $\sigma^2$ is the variance of the noise \cite{huber2011robust}. 

\begin{figure}
\centering
\includegraphics[width=3in,angle=0]{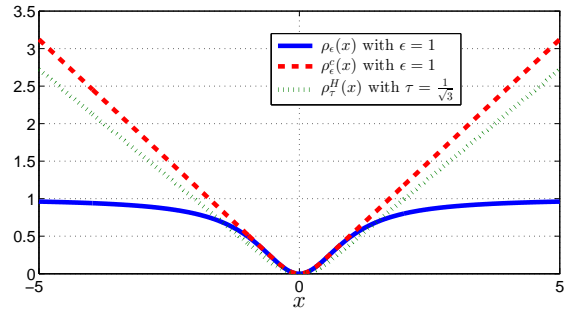}\vspace{-1.5mm}
\caption{\small{ Comparison of the IRLS weight function, its convex approximation, and the Huber norm. }}
\label{rho-fig}
\vspace{-5mm}
\end{figure}

The Huber norm is a convex function. To use the results of robust statistics in the proposed problem, a convex version of the cost function in (\ref{Rho}) should be employed. The function $\rho_\epsilon(.)$, can be surrogated by its closest convex approximation, 
\begin{equation}
\scriptsize
\label{rho-convex}
\rho_\epsilon^c(x) = \left\{ \begin{array}{lr} \frac{x^2}{x^2 + \epsilon^2}  & : | x | < \epsilon_0 \\[2ex] \frac{1}{8}(\frac{3}{\epsilon_0} |x| - 1 ) & : | x | \geq \epsilon_0 \end{array} \right. 
\end{equation}
with $\epsilon_0 = \frac{\epsilon}{\sqrt[]{3}}$. Figure \ref{rho-fig} illustrates the similarity between the Huber norm and the convex approximation of $\rho_\epsilon(.)$, i.e., $\rho_\epsilon^c(x)$. The cost functions resemble a least square estimator for errors less than a cut-off parameter, which is the optimal cost function for Gaussian noise. On the other hand, for large values of error, the cost functions resemble the $\ell_0$ or $\ell_1$ norms, which are known to promote sparsity.\par 

By extending the results of robust statistics to the proposed problem, we utilize the same cut-off parameter for $\rho_\epsilon^c(x)$ as the Huber norm. It means that for the case of Gaussian noise, we set $\epsilon = 1.34\; \sqrt[]{3} \; \sigma$, assuming that the nominal noise variance is available. If $\sigma$ is unknown, an estimation of it can be used \cite[Sec. 4.4]{maronna2006robust}. The numerical experiments in Section \ref{results} show that the estimator meets the Cr\`amer-Rao lower bound for sufficiently large number of sensors, by setting $\epsilon = 1.34\; \sqrt[]{3} \; \sigma$. \par 
\vspace{-2mm}
 \section{Proof of Theorem \ref{convergencetheorem}}\label{conv-analysis}
 \vspace{-1mm}
\label{conv-analysis}
Algorithm \ref{SRIRLSalg} alternates between two subproblems introduced in (\ref{Opt2}) and (\ref{SRL1matrix}). As discussed in Section \ref{SR-IRLS}, the optimization problem in (\ref{SRL1matrix}) is a GTRS and has a global minimizer for all the iterations. Moreover, $\boldsymbol{y}^{(k)}$, the global minimizer of (\ref{SRL1matrix}), is obtained by exploiting the conditions in (\ref{GTRS}). \par 

Also the optimization problem in (\ref{Opt2}) is strictly convex and the global minimizer, $\boldsymbol{w}^{(k)}$, can be calculated using the update rule in (\ref{SRIRLSW}) at each iteration. \par


\begin{lem}
\label{lem-nonincreasing}
$\{\mathcal{J}(\boldsymbol{y}^{(k)},\boldsymbol{w}^{(k)} ) \}$ is non-increasing using the update rules in Algorithm \ref{SRIRLSalg}, i.e., $\mathcal{J}(\boldsymbol{y}^{(k+1)},\boldsymbol{w}^{(k+1)} ) \leq \mathcal{J}(\boldsymbol{y}^{(k)},\boldsymbol{w}^{(k)} ), \forall k = 1,2,\ldots$
\end{lem}

\begin{IEEEproof}
Using the update rules in Algorithm \ref{SRIRLSalg}, we have
$$
\small
\mathcal{J}(\boldsymbol{y}^{(k+1)},\boldsymbol{w}^{(k+1)} ) \leq \mathcal{J}(\boldsymbol{y}^{(k+1)},\boldsymbol{w}^{(k)} ) \leq \mathcal{J}(\boldsymbol{y}^{(k)},\boldsymbol{w}^{(k)} ).
$$

The first inequality uses the fact that ${w}^{(k+1)}$ is the global minimizer of $\mathcal{J}(\boldsymbol{y}^{(k+1)},\boldsymbol{w} )$. Likewise, the second inequality uses the fact that ${y}^{(k+1)}$ is the global minimizer of $\mathcal{J}(\boldsymbol{y},\boldsymbol{w}^{(k)} )$. 
\end{IEEEproof}

Since $\mathcal{J}(\boldsymbol{y}^{(1)},\boldsymbol{w}^{(0)} ) < \infty$ and $\mathcal{J}(\boldsymbol{y}^{(k)},\boldsymbol{w}^{(k)} )$ is non-increasing, then either $\{ \mathcal{J}(\boldsymbol{y}^{(k)},\boldsymbol{w}^{(k)} ) \} \rightarrow -\infty$ , or $\{ \mathcal{J}(\boldsymbol{y}^{(k)},\boldsymbol{w}^{(k)} ) \}$ converges to some limit and  $\{ \mathcal{J}(\boldsymbol{y}^{(k+1)},\boldsymbol{w}^{(k+1)} ) - \mathcal{J}(\boldsymbol{y}^{(k)},\boldsymbol{w}^{(k)} ) \rightarrow 0 \}$ as $k \rightarrow \infty$.

Here, by setting the constant $\epsilon > 0$, we can assure that $-\ln w_i > -\infty, \forall i$. Then, it is easy to notice that $\{ \mathcal{J}(\boldsymbol{y}^{(k)},\boldsymbol{w}^{(k)} ) \}$ is bounded and the sequence $\{ \mathcal{J}(\boldsymbol{y}^{(k)},\boldsymbol{w}^{(k)} ) \} $ will converge to a constant value. To study the convergence of the iterates $\{ \boldsymbol{y}^{(k)},\boldsymbol{w}^{(k)} \}$, the definition of a limit point is presented \cite{Razaviyayn13BCDMM}.  \par 

\begin{define}
\label{limitpointdef}
\textit {
$\bar{x}$ is a limit point of $\{x^{(k)}\}$ if there exists a subsequence of $\{x^{(k)}\}$ that converges to $\bar{x}$. Note that every bounded sequence in $\mathbb{R}^n$
 has a limit point (or convergent subsequence).
 }
\end{define}



Now there exist a subsequence $\{ (\boldsymbol{y}^{(k_s)},\boldsymbol{w}^{(k_s)})\}$ that converges to a limit point $ (\boldsymbol{y}^{*},\boldsymbol{w}^{*})$. By plugging in $\boldsymbol{y}^{*}$ and $\boldsymbol{w}^{*}$ into the update rules, we will have 
$$
\boldsymbol{A}^T \boldsymbol{W}^{*} (\boldsymbol{A}\boldsymbol{y}^{*} - \boldsymbol{b}) + \lambda^{*} (\boldsymbol{D}\boldsymbol{y}^{*} + \boldsymbol{f}) =  0,  
$$
$$
(\boldsymbol{\tilde{a_i}}^T\boldsymbol{y}^{*} - b_i)^2 + \epsilon^2 - \frac{1}{w_{i}^{*}} = 0 , \forall i,
$$
which are the derivatives of the Lagrange function of (\ref{Eq-FullOpt}) w.r.t. $\boldsymbol{y}$ and $w_i$. Thus, $(\boldsymbol{y}^{*},\boldsymbol{w}^{*})$ is a stationary point of (\ref{Eq-FullOpt}). \par 
\scriptsize
{
\balance
\bibliography {Mendeley}
\bibliographystyle{ieeetr}
}

\end{document}